\documentclass{amsart}
\usepackage{graphicx}
\usepackage[all]{xy}
\usepackage{latexsym}
\usepackage{amssymb}
\usepackage{amsmath}
\usepackage{color}

\xyoption{arc}

 \newfam\cyrfam

  \font\tencyr=wncyr10

  \font\sevencyr=wncyr7

  \font\fivecyr=wncyr5

  \textfont\cyrfam=\tencyr \scriptfont\cyrfam=\sevencyr

    \scriptscriptfont\cyrfam=\fivecyr

  \newfam\cyifam

  \font\tencyi=wncyi10

  \font\sevencyi=wncyi7

  \font\fivecyi=wncyi5

  \textfont\cyifam=\tencyi \scriptfont\cyifam=\sevencyi

    \scriptscriptfont\cyifam=\fivecyi

\def\id{{\mbox{1 \hskip -7pt 1}}}
\newcommand{\sgn}{{\mathit s  \mathit g\mathit  n}}
 \newcommand{\lon}{\longrightarrow}
 \newcommand{\bu}{\bullet}
 
 \newcommand{\rar}{\rightarrow}

\newcommand{\tw}{\mathsf{tw}}
\newcommand{\RGra}{\mathcal{R} \mathcal{G} ra}
\newcommand{\GRav}{\mathcal{G} \mathcal{R}av}

 \newcommand{\Z}{{\mathbb Z}}
 \newcommand{\bS}{{\mathbb S}}

 \newcommand{\N}{{\mathbb N}}
 \newcommand{\K}{{\mathbb K}}

 \newcommand{\ot}{\otimes}
 
 \newcommand{\Beq}{\begin{equation}}
 \newcommand{\Eeq}{\end{equation}}
 \newcommand{\Beqr}{\begin{eqnarray}}
 \newcommand{\Eeqr}{\end{eqnarray}}
 \newcommand{\Beqrn}{\begin{eqnarray*}}
 \newcommand{\Eeqrn}{\end{eqnarray*}}
 \newcommand{\Ba}{\begin{array}}
 \newcommand{\Ea}{\end{array}}
 \newcommand{\Bi}{\begin{itemize}}
 \newcommand{\Ei}{\end{itemize}}
 \newcommand{\Bc}{\begin{center}}
 \newcommand{\Ec}{\end{center}}

 \newcommand{\f}{{\mathcal O}}

 \newcommand{\cE}{{\mathcal E}}
 
 \newcommand{\cG}{{\mathcal G}}
 \newcommand{\caH}{{\mathcal H}}

 \newcommand{\caL}{{\mathcal L}}
 \newcommand{\cM}{{\mathcal M}}

 \newcommand{\cR}{{\mathcal R}}
 
 \newcommand{\cT}{{\mathcal T}}

 \newcommand{\ga}{\gamma}
 
 \newcommand{\Ga}{\Gamma}

 \newcommand{\Hom}{{\mathrm H\mathrm o\mathrm m}}
 
 \newcommand{\sip}{\smallskip}
 \newcommand{\bip}{\bigskip}
 \newcommand{\mip}{\vspace{2.5mm}}

\newcommand{\LB}{\mathcal{L}\mathit{ieb}}

\newcommand{\Holie}{\mathcal{H} \mathit{olie}}
\newcommand{\Lie}{\mathcal{L} \mathit{ie}}

\newcommand{\DRGra}{\mathcal{RG}ra^\uparrow}
\newcommand{\RGraphs}{\mathcal{RG}raphs}

\newcommand{\ORGra}{\mathcal{RG}ra^{or}}
\newcommand{\ORGraphs}{\mathcal{RG}raphs^{or}}
\newcommand{\OORGraphs}{O\mathcal{RG}raphs}
\newcommand{\DRTrees}{\mathcal{RT}rees^\uparrow}
\newcommand{\DRTre}{\mathcal{RT}r^\uparrow}

\theoremstyle{plain}
\swapnumbers

\newtheorem{prop-def}[theorem]{Proposition-definition}

\newtheorem{f-theorem}{Formality Theorem}[section]
\newtheorem{main-theorem}{Main~Theorem}[section]
\newtheorem{section-theorem}{Theorem}[section]

\theoremstyle{definition}

\usepackage{tikz}
\usetikzlibrary{matrix}
\usetikzlibrary{shapes}
\usetikzlibrary{arrows}
\usetikzlibrary{calc,3d}
\usetikzlibrary{decorations,decorations.pathmorphing,decorations.pathreplacing,decorations.markings}
\usetikzlibrary{through}

\tikzset{ext/.style={circle, draw,inner sep=1pt},int/.style={circle,draw,fill,inner sep=1.4pt},nil/.style={inner sep=1pt}}
\tikzset{cy/.style={circle,draw,fill,inner sep=2pt},scy/.style={circle,draw,inner sep=2pt},scyx/.style={draw,cross out,inner sep=2pt},scyt/.style={draw,regular polygon,regular polygon sides=3,inner sep=0.95pt}}
\tikzset{exte/.style={circle, draw,inner sep=3pt},inte/.style={circle,draw,fill,inner sep=3pt}}
\tikzset{diagram/.style={matrix of math nodes, row sep=3em, column sep=2.5em, text height=1.5ex, text depth=0.25ex}}
\tikzset{diagram2/.style={matrix of math nodes, row sep=0.5em, column sep=0.5em, text height=1.5ex, text depth=0.25ex}}

 \tikzset{
  rightblack/.style={
    decoration={markings,mark=at position .8 with {\arrow[scale=1.2,black]{latex}}},
    postaction={decorate},
    shorten >=0.4pt}}
\tikzset{
  leftblack/.style={
    decoration={markings,mark=at position .55 with {\arrowreversed[scale=1.2,black]{latex}}},
    postaction={decorate},
    shorten >=0.4pt}}

\begin{document}

 \sloppy

 \newenvironment{proo}{\begin{trivlist} \item{\sc {Proof.}}}
  {\hfill $\square$ \end{trivlist}}

\long\def\symbolfootnote[#1]#2{\begingroup%
\def\thefootnote{\fnsymbol{footnote}}\footnote[#1]{#2}\endgroup}

\title{Pre-Calabi-Yau algebras and oriented gravity properad}

\author{Sergei Merkulov\\
{$\text{\hspace{90mm}\it To the memory of Jan-Erik Roos}$}
}
\address{Sergei~A.\ Merkulov:
Department of Mathematics, Luxembourg University, Grand Duchy of Luxembourg}
\email{smerkulov25@gmail.com}

\begin{abstract}  We study the dual cyclic Hochschild  complex $Cyc^\bu(A,\K)$ of a
(possibly, infinite-dimensional) $A_\infty$-algebra $(A,\mu)$ and prove
that any pre-Calabi-Yau extension  $\pi$ of the given $A_\infty$ structure $\mu$ in $A$
induces on the cyclic cohomology of $(A,\mu)$ a representation of a new dg properad of {\em oriented}\, ribbon graphs. We compute the cohomology of that properad in terms of the compactly supported cohomology groups of moduli spaces
$\cM_{g,m+n}$ of algebraic curves of genus $g$ with $m\geq 1$ boundaries and   $n\geq 1$
marked points.
We also show that the gravity operad acts naturally
on the higher Hochschild cohomology of any pre-CY algebra $(A, \pi)$.


\end{abstract}

 \maketitle
\markboth{}{}

{\Large
\section{\bf Introduction}
}

The authors of \cite{KTV} have introduced a new dg prop of certain ribbon quivers $Q^d=\{Q^d(m,n)\}_{m,n\geq 1}$ which has the following main properties:
\Bi
\item[(i)] For any pre-Calabi-Yau extension of a given $A_\infty$ structure $\mu$ on a graded
vector space $A$, there is an associated  natural action of $Q^d$ on the Hochschild
chain complex
$$
C_\bu(A) =\prod_{k\geq 1}  (A[1])^{\ot k}
$$
of $(A,\mu)$.

\item[(ii)] The cohomology groups $H^\bu(Q^d(m,n))$ of the dg prop $Q^d$ are isomorphic to the totality  $\prod _{g\geq 0} H^\bu(\cM_{g, \vec{m}+\vec{n}}, \caL^d)$  of the cohomology groups (with coefficients in powers of a certain rank one local system $\caL$) of the moduli spaces $\cM_{g, \vec{m}+\vec{n}}$ of genus $g$ algebraic curves
with two sets of labelled marked points, one of cardinality $m\geq 1$ and another of cardinality $n\geq 1$, which are equipped with tangent rays.

\Ei

\sip

In this note we show a version of the above results for the dual {\em cyclic}\,  Hochschild complex of $(A,\mu)$,
$$
Cyc^\bu(A,\K) :=\prod_{k\geq 1} \Hom\left( (A[1])^{\ot k}, \K\right)_{\Z_k}.
$$

Our starting point is a properad of connected ribbon graphs $\RGra_d$ introduced in \cite{MW1} which comes equipped with a canonical map
$$
\LB_d \lon \RGra_d
$$
from the properad $\LB_d$ of (degree shifted) Lie bialgebras. In particular, $\RGra_d$ is a properad
under the operad of (degree shifted) Lie algebras $\Lie_d$ and hence can be twisted into a
{\em differential}\, graded properad $(\tw \RGra_d, \delta)$ via a general construction
introduced in \cite{Wi1}. Let
$$
\RGraphs_d=\{\RGraphs_d(m,n), \delta\}_{m\geq1, n\geq 0}
$$
be the sub-properad of $\tw \RGra_d$ generated by ribbon graphs whose unlabelelled vertices are at
least trivalent; the number $n$ stands for the number of {\em labelled}\, vertices of the generating ribbon graphs, while $m$ for the number of their (all labelled) boundaries; $\delta$ is the differential which acts on ribbon graphs by splitting labelled and unlabelled vertices (its explicit formula is given in (2) below). The complex $(\RGraphs_d(m,0), \delta)$ is precisely the (degree shifted) classical Penner's complex
 \cite{Pe,Ko1} of ribbon graphs with marked boundaries,
so that one has
$$
H^\bu(\RGraphs_d(m,0), \delta)\  =\prod_{g\geq 0\atop 2g+m\geq 3}
 H^{\bu-m+d(2g-2+m)}_c(\cM_{g,m})
$$
where $\cM_{g,m}$ is the moduli spaces of genus $g$ algebraic curves with $m$ marked points,
$H^\bu_c$ stands for the compactly supported cohomology functor.
More generally,  for any $m\geq 1$, $n\geq 0$ one has \cite{Co,Me1}
$$
H^\bu(\RGraphs_d(m,n), \delta) = \prod_{g\geq 0\atop 2g+m+n\geq 3}
H_c^{\bu -m +d(2g-2+m+n) }(\cM_{g, m+n}),
$$
 where $\cM_{g,m+n}$ be the moduli space of algebraic curves of genus $g$ with $m\geq 1$ boundaries
 and   $n\geq 0$ marked points (or, equivalently, the moduli space of algebraic curves of genus
 $g$ with $m+n\geq 1$ marked points,
  $m\geq 1$ of them are called out-points, and the remaining $n\geq 0$ are called in-points).
  Hence the  above collection of cohomology groups of moduli spaces $\cM_{g,m+n}$ has the structure
  of a properad which is called in \cite{Me1} the {\em gravity properad} as its
genus zero part coincides precisely with Getzler's gravity operad \cite{Ge}.

\sip

Using a version of the properad $\RGra_{d}$ with {\em directed}\, edges we construct in \S 3 of
this paper a new dg properad of {\em oriented}\, ribbon graphs (or, equivalently, of
ribbon {\em quivers}),
$$
\OORGraphs_d=\{\OORGraphs_d(m,n), \delta\}_{m\geq1, n\geq 1},
$$
and prove the following two main theorems of this paper.

\sip

\subsection{Theorem}\label{1: Main theorem on action of OOGraphs on Cyc} {\em Let $(A,\mu)$ be a (possibly infinite-dimensional) $A_\infty$-algebra. Then any degree $d$ pre-Calabi-Yau extension $\pi$  of the given $A_\infty$-algebra structure $\mu$ in $A$ induces a natural action of the dg properad $\OORGraphs_d$ on the dual cyclic Hochschild complex $Cyc^\bu(A,\K)$ of $(A,\mu)$.
}

\sip

The dg properad $\OORGraphs_{d+1}$ comes equipped with a morphism from the minimal resolution
$\Holie_d$
of the operad of (degree shifted) $\Lie_d$ algebras,
$$
\Holie_d \lon \OORGraphs_{d+1},
$$
which is non-trivial on every generator of $\Holie_d$
(see Proposition {\ref{3: prop on F from holieb to RGraphs}} for an explicit formula).

\subsection{Corollary}\label{1: corollary on Holie} {\em Any degree
$(d+1)$ pre-Calabi-Yau extension $\pi$ of a given $A_\infty$-algebra structure $\mu$ in
$A$ makes the dual cyclic Hochschild complex $Cyc^\bu(A,\K)$ into a $\Holie_d$-algebra.
}

\sip

The main result of this paper is a purely combinatorial proof (see \S 4 below) of the following claim.

\subsection{Main Theorem}\label{1: Main theorem on H(OOGraphs)} {\em For any $d\in \Z$
and any $m,n\geq 1$ one has}
$$
H^\bu(\OORGraphs_{d+1}(m,n))=H^\bu(\RGraphs_d(m,n))=
\prod_{g\geq 0\atop 2g+m+n\geq 3}
H_c^{\bu -m +d(2g-2+m+n) }(\cM_{g, m+n})
$$
Notice the degree shift $d+1$ in the first cohomology group. Due to this result, one may call
$\OORGraphs_{d+1}$ the {\em oriented}\, gravity properad; it gives us a new combinatorial model
for the totality of the cohomology groups $H^\bu_c(\cM_{g, m+n})$.

\sip

We also show that the cohomology $H^\bu(C_{[d]}^\bu(A), \delta_\pi)$ of the higher Hochschild complex of any pre-CY algebra $(A,\pi)$ is naturally an algebra over the gravity operad, with the binary operation being the one induced from the standard Lie bracket
 in $C_{[d]}^\bu(A)$. A gravity algebra structure on the negative cyclic homology of
Calabi-Yau algebras was established earlier in \cite{CEL}.

  \subsection{Some notation} We work over a field $\K$ of characteristic zero.
 The set $\{1,2, \ldots, n\}$ is abbreviated to $[n]$;  its group of automorphisms is
denoted by $\bS_n$; the trivial (resp., the sign) one-dimensional representation of
 $\bS_n$ is denoted by $\id_n$ (resp.,  $\sgn_n$). The cardinality of a finite set $S$ is denoted by $\# S$ while its linear span over a
field $\K$ by $\K\left\langle S\right\rangle$.


\sip

If $V=\oplus_{i\in \Z} V^i$ is a graded vector space, then
$V[k]$ stands for the graded vector space with $V[k]^i:=V^{i+k}$. For $v\in V^i$ we set $|v|:=i$.

\sip

For a
(pr)operad $\f$ we denote by $\f\{k\}$ a (pr)operad which is uniquely defined by
 the following property:
for any graded vector space $V$ a representation
of $\f\{k\}$ in $V$ is identical to a representation of  $\f$ in $V[k]$.
 The degree shifted operad of Lie algebras $\caL \mathit{ie}\{d\}$  is denoted by $\caL ie_{d+1}$ while its minimal resolution by $\caH \mathit{olie}_{d+1}$; representations of $\caL ie_{d+1}$ are vector spaces $V$ equipped with the Lie bracket of degree $-d$.

\sip

{\bf Acknowledgement}. It is a great pleasure to thank Thomas Willwacher for many valuable
discussions.

\bip

{\Large
\section{\bf Reminder on pre-Calabi-Yau algebras}
}

\sip

\subsection{A list of main graded vector spaces under study } Let $A$ be a (possibly infinite-dimensional) graded vector space over a
field $\K$. We consider
\Bi
\item[(i)]
$$
C^\bu(A,A):= \prod_{k\geq 1} \Hom\left( (A[1])^{\ot k}, A[1]\right),
$$
\item[(ii)]
$$
C^\bu(A,\K) :=\prod_{k\geq 1} \Hom\left( (A[1])^{\ot k}, \K\right)
$$
\item[(iii)]
$$
Cyc^\bu(A,\K) :=\prod_{k\geq 1} \Hom\left( (A[1])^{\ot k}, \K\right)_{\Z_k}
$$
the graded vector space of (co)invariants under the natural actions of
the cyclic groups $\Z_k$ on $\Hom\left( (A[1])^{\ot k}, \K\right)$.
\item[(iv)]
$$
C_{d}^\bu(A,A) :=\prod_{k\geq 1} C_d^{(k)}(A,A):= \prod_{k\geq 1}
\left( \bigoplus_{n_1,\ldots, n_k\geq 0}
\Hom\left(\bigotimes_{i=1}^{k} (A[1])^{ \ot n_i}, (A[2-d])^{\ot k}\right)\right)
$$
\item[(v)]
$$
C_{[d]}^\bu(A) :=\prod_{k\geq 1} C_{[d]}^{(k)}(A):= \prod_{k\geq 1}
\left( \bigoplus_{n_1,\ldots, n_k\geq 0}
\Hom\left(\bigotimes_{i=1}^{n_i}
(A[1])^{\ot n_i}, (A[2-d])^{\ot k}\right)_{\Z_k}\right)
$$
the graded vector space (called the {\em higher}\,
 Hochschild space of $A$) of (co)invariants under the natural actions of the cyclic
 groups $\Z_k$ on
$\Hom\left(\bigotimes_{i=1}^{k} (A[1])^{\ot n_i}, (A[2-d])^{\ot k}\right)$.

\item[(vi)]
$$
\widehat{C}_{[d]}^\bu(A) :=
Cyc^\bu(A,\K) \oplus  C_{[d]}^\bu(A)
$$
called the {\em full}\, higher  Hochschild space of $A$.  If $A$ is finite-dimensional, then
$\widehat{C}_{[d]}^\bu(A)$ can be identified with  the graded vector space
$$
 Cyc(W):= \prod_{k\geq 1} (\ot^k W)_{Z_k}, \ \ \
  \text{where}\ W:= A[2-d]\oplus A^*[-1]
$$
while
$C_{[d]}^\bu(A)$ gets identified with the subspace of the latter which is
generated by cyclic words having at least one letter from the direct summand $A[2-d]\subset W$.
Note that $W$ comes equipped with a natural non-degenerate pairing in degree $1-d$.

\Ei

\subsection{Remark on degree shift conventions} The above definitions agree with the ones in \cite{KTV} if one replaces the integer parameter $d$ to $4-d$  and then shifts the degree $C^\bu_{[d]}(A) \rar C^\bu_{[d]}[2-d]$.

\sip
Our choice is consistent with the well-established
role of the parameter $d$ in the theory of graph complexes: the natural composition of two maps from, say, $C_{[d]}^{(1)}(A,A)$, has degree $1-d$.

\subsection{Higher Hochschild complexes of an $A_\infty$ algebra}
 Let $A$ be a (possibly infinite-dimensional) $A_\infty$ algebra, that is, $A$ comes equipped
 with a collection of linear maps
$$
\mu=\left\{\mu_n: \ot^n (A[1]) \lon A[2]\right\}_{n\geq 1}
$$
satisfying well-known equations.  This structure makes  all the graded vector spaces in (i)-(vi) into
{\em complexes}\, equipped with the Hochschild differential which we denote by $d_\mu$.

\sip

\subsection{Pre-CY algebras}

 The graded vector space $C_{[d]}^\bu(A)$ has a natural
$\Lie_{d}$-algebra structure \cite{IK,IKV}. Its Maurer-Cartan elements $\pi$, that is
 degree $d$ elements $\pi\in C_{[d]}^\bu(A)$ satisfying the equation
$$
[\pi,\pi]=0,
$$
are called  {\em pre-CY structures on $A$ of degree $d$} in \cite{IK,IKV, KTV}.
If $A$ is finite-dimensional, such MC elements have been studied in \cite{TZ} under the name
  {\em $V_\infty$-algebras}. We refer to
\S 1.4 in \cite{KTV} for more details about the history of this notion.

\sip

If $\pi$ is a pre-CY algebra structure on $A$, then the summand $\mu$  of $\pi$
which belongs to the direct summand $C_{[d]}^{(1)}\subset C_{[d]}^\bu(A)$ defines an
$A_\infty$-structure on $A$, and then  $\pi$ is called a {\em pre-CY extension}\,
of $\mu$. If $(A,\pi)$ is a pre-CY algebra, then the higher Hochschild complex $C_{[d]}^\bu(A)$
has two natural differentials, $\delta_\mu=[\mu, ...]$ and  $\delta_\pi=[\pi,...]$. The dg Lie algebras
$$
\left({C}_{[d]}^\bu(A), [\ ,\ ],\delta_\mu\right) \ \ \& \ \
\left({C}_{[d]}^\bu(A), [\ ,\ ],\delta_\pi\right).
$$
are called the {\em higher Hochschild complexes}\, of $\mu$ and, respectively, of $\pi$. We study below
natural  operations (controlled by certain properads of ribbon graphs) on both these complexes
which are associated with any pre-CY extension $\pi$ of the given $A_\infty$-structure $\mu$.

\sip
A nice dg properad controlling an important class of pre-Calabi-Yau structures (the ones which can be called {\it strongly homotopy}\,  double Poisson structures)  has been explicitly described and  studied 
in \cite{LV}.

\subsection{Curved pre-CY algebras} The graded vector space $\widehat{C}_{[d]}^\bu(A)$ is also
naturally a $\Lie_{d}$-algebra; its Maurer-Cartan elements $\hat{\pi}\in \widehat{C}_{[d]}^\bu(A)$ 
can be called  {\em extended}\, or {\em  curved}\, pre-CY structures in $A$. Such structures are not studied
   in this paper, but sometimes we use the underlying (so called) {\em full higher Hochschild complexes},
$$
\left(\widehat{C}_{[d]}^\bu(A), [\ ,\ ],\delta_\mu\right) \ \ \& \ \
\left(\widehat{C}_{[d]}^\bu(A), [\ ,\ ],\delta_\pi\right),
$$
associated with a {\em non-curved}\, $A_\infty$- structure $\mu$ in $A$ and its
{\em non-curved}\, pre-CY extension $\pi$.

\bip


{\Large
\section{\bf Properad of oriented ribbon graphs and higher Hochschild complexes}
}

\mip

\subsection{Reminder on the properad  of (undirected) ribbon graphs}
Let $\RGra_d=\{\RGra_d(m,n)\}_{m,n\geq 1}$ be a properad of {\em connected}\, ribbon graphs
introduced\footnote{More precisely, the symbol $\RGra_d$ stands in \cite{MW1} for the {\em prop}\, 
generated by not necessarily connected ribbon graphs; in this paper we work solely with connected
graphs and hence use the symbol $\RGra_d$ for the sub-properad of the latter which is generated 
by connected ribbon graphs.} in \S 4 of \cite{MW1}\footnote{
It is worth noting that there are other constructions of prop(erad)s based 
on the idea of ribbon graph
 which can be found in the literature (see e.g. \cite{KS, KTV,TZ,WW}); the properad $\RGra_d$ 
 and its twisted version $\tw\RGra_d$ studied in \cite{MW1,Me1} are different  from that 
 constructions because the automorphism groups and orientations of the generating ribbon graphs
  are quite  different. The aforementioned properads of ribbon graphs
   control different algebraic structures.}. The $\bS_m\times\bS_n$-module $\RGra_d(m,n)$
is generated by ribbon graphs $\Ga$ with $n$ labelled vertices and $m$ labelled
(by integers $\bar{1}, \bar{2},\ldots$) boundaries, e.g.
$$
\Ba{c} \resizebox{11mm}{!}{\xy
 (5,3)*{^{\bar{1}}};
 (10,1)*+{_1}*\frm{o}="B";
 (0,1)*+{_2}*\frm{o}="A";
 \ar @{.>} "A";"B" <0pt>
\endxy} \Ea
\in \RGra_d(1,2), \ \
\Ba{c}\resizebox{17mm}{!}{
\mbox{$\xy
 (-2,3)*{^{\bar{1}}};
 (-2,-3)*{_{\bar{2}}};
 (-9,0)*{^{\bar{3}}};
 (-4,0)*+{_{_4}}*\frm{o}="C";
  (9,0)*+{_{_1}}*\frm{o}="1";
(-7,8)*+{_{_2}}*\frm{o}="2";
(-7,-8)*+{_{_3}}*\frm{o}="3";
 \ar @{.>} "1";"C" <0pt>
  \ar @{.>} "1";"2" <0pt>
   \ar @{.>} "1";"3" <0pt>
 \ar @{.>} "2";"C" <0pt>
  \ar @{.>} "3";"C" <0pt>
\endxy$}}
\Ea
\in \RGra_d(3,4), \ \
\Ba{c}\resizebox{12mm}{!}{\xy
(0,0)*{_{\bar{1}}};
 (0,5)*{_{\bar{2}}};
   {\ar@{.>}@/^0.7pc/(-6,0)*+{_1}*\frm{o} ;(6,0)*+{_2}*\frm{o}};
 {\ar@{<.}@/^0.7pc/(6,0)*+{_2}*\frm{o};(-6,0)*+{_1}*\frm{o}};
\endxy}\Ea\in \RGra_d(2,2).
$$
Each edge of $\Ga$ is equipped with a direction which can be flipped with the
 sign factor,
\Beq\label{2: symmetry of dotted edges}
\Ba{c} \resizebox{11mm}{!}{\xy
 (10,1)*+{_k}*\frm{o}="B";
 (0,1)*+{_i}*\frm{o}="A";
 \ar @{.>} "A";"B" <0pt>
\endxy} \Ea
= (-1)^{d}
\Ba{c} \resizebox{11mm}{!}{\xy
 (10,1)*+{_k}*\frm{o}="B";
 (0,1)*+{_i}*\frm{o}="A";
 \ar @{<.} "A";"B" <0pt>
\endxy} \Ea.
\Eeq
Hence we skip from now on showing directions on dotted edges (assuming tacitly that some 
choice has been made). The cohomological degree of $\Ga$ is defined by
$$
|\Ga|=(1-d) \# E(\Ga)
$$
where $E(\Ga)$ stands for the set of edges, i.e.\ each edges is assigned the degree $1-d$;
for $d$ even it is assumed that some ordering of edges of $\Ga$ is fixed up to an 
even permutation (an odd permutation acts as the multiplication by $-1$).

\sip

The properadic compositions in $\RGra_d$ are given by substituting a boundary $b$ of one
ribbon graph into a vertex $v$ of another ribbon graph and reattaching edges glued 
earlier to $v$ among
(the arcs of) the vertices belonging to the boundary $b$ in all possible ways while 
respecting the obvious
cyclic orderings of edges at $v$ and (the arcs of the) vertices in $b$ (see \S 4.2 in 
\cite{MW1} for full details and illustrating examples).

\sip

Let $W$ be an arbitrary graded vector space and  $
Cyc^\bu W=\oplus_{n\geq 1} (W^{\ot n})^{\Z_n}
$ the associated space of cyclic words. Then any  (not necessarily non-degenerate)
 pairing $\Theta:  W\ot W  \lon  \K[1-d]$ satisfying the symmetry condition,
$$
\Theta(w_1,w_2)=(-1)^{d+|w_1||w_2|}\Theta(w_2,w_1).
$$
gives rise to a morphism of properads,
$$
\rho_W: \cR \cG ra_{d} \lon \cE nd_{Cyc^\bu W}.
$$
i.e.\ defines a representation of $\RGra_d$ in $Cyc^\bu(W)$ (see \S 4.2.2 in \cite{MW1} 
for full details).

\sip

If $\Lie_d$ stands for the operad of degree shifted Lie algebras,
then there is a  morphism of properads \cite{MW1}
$$
i: \Lie_d \lon    \RGra_d
$$
given on the Lie bracket generator of $\Lie_d$  by
$$
\Ba{c}\begin{xy}
 <0mm,0.66mm>*{};<0mm,3mm>*{}**@{-},
 <0.39mm,-0.39mm>*{};<2.2mm,-2.2mm>*{}**@{-},
 <-0.35mm,-0.35mm>*{};<-2.2mm,-2.2mm>*{}**@{-},
 <0mm,0mm>*{\bu};<0mm,4.1mm>*{^{^{\bar{1}}}}**@{},
   <0.39mm,-0.39mm>*{};<2.9mm,-4mm>*{^{_2}}**@{},
   <-0.35mm,-0.35mm>*{};<-2.8mm,-4mm>*{^{_1}}**@{},
\end{xy}\Ea
 \stackrel{i}{\lon}
 \Ba{c} \resizebox{11mm}{!}{\xy
 (10,1)*+{_2}*\frm{o}="B";
 (0,1)*+{_1}*\frm{o}="A";
 \ar @{.} "A";"B" <0pt>
\endxy} \Ea
$$
Hence one can apply to $\RGra_d$ Thomas Willwacher's  \cite{Wi1}
twisting endofunctor $\tw$  to obtain\footnote{In this particular context, see also 
an overview of the endofunctor $\tw$ in terms of decorated
graphs given in \cite{Me2}.} a dg properad
$$
\tw \RGra_d=\left\{ \tw \RGra_d(m,n), \delta   \right\}_{m\geq 1, n \geq 0}
$$
which is generated by ribbon graphs $\Ga$ with $m\geq 1$ labelled boundaries, $n\geq 0$ 
labelled vertices, and any
number of unlabelled vertices (which we denote in pictures as small white circles 
$\circ$) to which we assign the cohomological degree $d$ (see \S 3.3 in \cite{Me1} 
for more details). For example,
$$
\Ba{c} \resizebox{11mm}{!}{\xy
 (5,3)*{^{\bar{1}}};
 (10,1)*+{_1}*\frm{o}="B";
 (0,1)*+{_2}*\frm{o}="A";
 \ar @{.} "A";"B" <0pt>
\endxy} \Ea
\in \tw \RGra_d(1,2), \ \
\Ba{c}\resizebox{17mm}{!}{
\mbox{$\xy
 (-2,3)*{^{\bar{1}}};
 (-2,-3)*{_{\bar{2}}};
 (-9,0)*{^{\bar{3}}};
 (-4,0)*{\circ}="C";
  (9,0)*+{_{_1}}*\frm{o}="1";
(-7,8)*+{_{_2}}*\frm{o}="2";
(-7,-8)*{\circ}="3";
 \ar @{.} "1";"C" <0pt>
  \ar @{.} "1";"2" <0pt>
   \ar @{.} "1";"3" <0pt>
 \ar @{.} "2";"C" <0pt>
  \ar @{.} "3";"C" <0pt>
\endxy$}}
\Ea
\in \tw \RGra_d(3,2), \ \
\Ba{c}\resizebox{12mm}{!}{\xy
(0,0)*{_{\bar{1}}};
 (0,5)*{_{\bar{2}}};
   {\ar@{.}@/^0.7pc/(-6,0)*+{_1}*\frm{o} ;(6,0)*{\circ}};
 {\ar@{.}@/^0.7pc/(6,0)*{\circ};(-6,0)*+{_1}*\frm{o}};
\endxy}\Ea\in \tw \RGra_d(2,1).
$$
The cohomological degree of $\Ga\in \tw\RGra_d$ is given by
$$
|\Ga|=(1-d) \# E(\Ga) + d\# V_\circ(\Ga),
$$
where $V_\circ(\Ga)$ stands for the set of unlabelled vertices.
The differential in $\tw\RGra_d(m,n)$ is given by the standard ``splitting of 
vertices" formula,
\Beq\label{2: delta in Tw(RGra)}
\delta\Ga:=
\sum_{i=1}^m \Ba{c}\resizebox{4mm}{!}{  \xy
 (0,7)*{\circ}="A";
 (0,0)*+{_1}*\frm{o}="B";
 \ar @{.} "A";"B" <0pt>
\endxy} \Ea
\  _1\hspace{-0.7mm}\circ_i \Ga\ \
- \ \ (-1)^{|\Ga|} 
\sum_{j=1}^n\Ga\  _j\hspace{-0.7mm}\circ_1
\Ba{c}\resizebox{4mm}{!}{  \xy
 (0,6)*{\circ}="A";
 (0,0)*+{_1}*\frm{o}="B";
 \ar @{.} "A";"B" <0pt>
\endxy} \Ea
\  -(-1)^{|\Ga|}\ \frac{1}{2} \sum_{v\in V_\circ(\Ga)} \Ga\circ_v  \left(\xy
 (0,0)*{\circ}="a",
(5,0)*{\circ}="b",
\ar @{.} "a";"b" <0pt>
\endxy\right)
\Eeq
where the symbol $_i\circ_j$ stands for the properadic composition of the 
$i$-labelled vertex of one graph with the $j$-labelled boundary of another graph 
(see \S 4.2.2 in \cite{MW1} for details and illustrating examples), and  the 
symbol $\Ga\circ_v  \left(\xy
 (0,0)*{\circ}="a",
(5,0)*{\circ}="b",
\ar @{.} "a";"b" <0pt>
\endxy\right)$ means the substitution of the graph  $\xy
 (0,0)*{\circ}="a",
(5,0)*{\circ}="b",
\ar @{.} "a";"b" <0pt>
\endxy$ into the unlabelled vertex $v$ of the graph $\Ga$ followed by
 the summation over all possible re-attachments of the half-edges attached earlier
 to $v$ among the two unlabelled vertices in a way which respects their cyclic
 ordering. Note that for almost all graphs the new univalent unlabelled vertices
 arising in the first summand of $\delta$ cancel out the new univalent
 unlabelled vertices arising in the second and the third part of that differential.

\sip

 The dg properad  $\tw\cR\cG ra_d$ contains a dg sub-properad $\RGraphs_d$
  spanned by ribbon graphs  with every unlabelled vertex having valency  $\geq 3$. 
  It was shown in \S 3.4 of
  \cite{Me1} that
 the inclusion
$$
\RGraphs_d(m,n) \lon \tw \RGra_d(m,n)
$$
is  a quasi-isomorphism for any $m,n\geq 1$, while  for $n=0$ on has
$$
H^\bu(\tw \RGra_d(m,0))= H^\bu\left(\RGraphs_d(m,0)\right) \oplus \
\bigoplus_{p\geq 1\atop
p\equiv 2d+1 \bmod 4} \K[-p(d-1)]
$$
where the summand $\K[-p(d-1)]$ is generated by the polytope-like ribbon graph 
with $p$ edges and $p$ bivalent unlabelled vertices.

\sip

The complex $(\RGraphs_d(m,0), \delta)$ generated by ribbon graphs with {\em all}\, vertices 
unlabelled is precisely the classical
(degree shifted) R.\ Penner's ribbon graph complex \cite{Pe} (see also \cite{Ko1})
 with marked boundaries,
so that its cohomology has a nice geometric meaning,
$$
H^\bu(\RGraphs_d(m,0), \delta)\  =\prod_{2g+m\geq 3} H^{\bu-m+d(2g-2+m)}_c(\cM_{g,m})
$$
where $\cM_{g,m}$ is the moduli spaces of genus $g$ algebraic curves with $m$ marked points, 
and $H^\bu_c$ stands for the compactly supported cohomology functor. Using K.\ Costello's
construction of  moduli spaces of nodal disks \cite{Co}
it is not hard to conclude (see \cite{Me1} for full details) that for any $m\geq 1$, 
$n\geq 0$ one has
$$
H^\bu(\RGraphs_d(m,n), \delta) =\prod_{2g+m+n\geq 3\atop  m\geq 1, n\geq 0}
H_c^{\bu -m +d(2g-2+m+n) }(\cM_{g, m+n})
$$
 where $\cM_{g,m+n}$ be the moduli space of algebraic curves of genus $g$ with $m\geq 1$
  boundaries and   $n\geq 0$ marked points (or, equivalently, the moduli space of 
  algebraic curves of genus $g$ with $m+n\geq 1$ marked points,
  $m\geq 1$ of them are called out-points, and the remaining $n\geq 0$ are called in-points).

  \sip

  The simplest possible moduli space of genus zero curves with three marked points 
  (which is just a point) contributes into the cohomology properad $H^\bu(\RGraphs_d)$
   in three different incarnations, as $H^\bu(\cM_{0, 1+2})=\K[d-1]$, as 
   $H_c^\bu(\cM_{0, 2+1})=\sgn_2[d-2]$ and as  $H_c^\bu(\cM_{0,3+0})=\sgn_3[d-3]$;
    its was proven in \cite{Me1} that these three classes generate a  
 sub-properad of quasi-Lie bialgebras. Moreover, the cohomology of the genus zero part 
 $\cR\cT rees_d$ of $\RGraphs_d$ is an operad isomorphic \cite{Wa} to
   the so called {\em gravity operad}\, introduced in \cite{Ge}. Hence the cohomology 
properad $H^\bu(\RGraphs_d)$ was called in \cite{Me1} the {\em gravity properad}\, 
 (and denoted there by $\GRav_d$).

\subsection{A dg properad  of  oriented ribbon graphs}\label{2: subsec on RGra^or and RGraphs^or}
Let $\DRGra_{d}=\{\DRGra_{d}(m,n)\}_{m,n\geq 1}$
be a directed  version of the properad $\RGra_{d}$ discussed above.
More precisely, the $\bS_m^{op}\times \bS_n$-module $\DRGra_{d}(m,n)$ is generated by
 ribbon graphs
with $n$ labelled vertices and $m$ labelled boundaries, and with {\em directed}\, edges
which are assigned the cohomological degree $1-d$.  In the pictures we show edges 
with a fixed direction as {\em solid arrows}, e.g.
 $$
\Ba{c} \resizebox{11mm}{!}{\xy
 (5,3)*{^{\bar{1}}};
 (10,1)*+{_1}*\frm{o}="B";
 (0,1)*+{_2}*\frm{o}="A";
 \ar @{->} "A";"B" <0pt>
\endxy} \Ea
\in \DRGra_d(1,2), \ \
\Ba{c}\resizebox{17mm}{!}{
\mbox{$\xy
 (-2,3)*{^{\bar{1}}};
 (-2,-3)*{_{\bar{2}}};
 (-9,0)*{^{\bar{3}}};
 (-4,0)*+{_{_4}}*\frm{o}="C";
  (9,0)*+{_{_1}}*\frm{o}="1";
(-7,8)*+{_{_2}}*\frm{o}="2";
(-7,-8)*+{_{_3}}*\frm{o}="3";
 \ar @{->} "1";"C" <0pt>
  \ar @{->} "1";"2" <0pt>
   \ar @{->} "1";"3" <0pt>
 \ar @{->} "2";"C" <0pt>
  \ar @{>} "3";"C" <0pt>
\endxy$}}
\Ea
\in \DRGra_d(3,4), \ \
\Ba{c}\resizebox{12mm}{!}{\xy
(0,0)*{_{\bar{1}}};
 (0,5)*{_{\bar{2}}};
   {\ar@{->}@/^0.7pc/(-6,0)*+{_1}*\frm{o} ;(6,0)*+{_2}*\frm{o}};
 {\ar@{->}@/^0.7pc/(6,0)*+{_2}*\frm{o};(-6,0)*+{_1}*\frm{o}};
\endxy}\Ea\in \DRGra_d(2,2).
$$
This properad contains a sub-properad
$\ORGra_d$ generated by ribbon graphs with no closed paths of directed edges;
for example, all the above ribbon graphs except the last one belong to $\ORGra_{d}$. The properad
$\ORGra_{d}$ contains in turn an sub-properad $\DRTre_d\subset \ORGra_d$
 generated by directed  connected ribbon graphs of genus zero, that is, by directed 
 ribbon trees; as ribbon trees have only one boundary, $\DRTre_d$ is an example of 
 an operad.

\subsubsection{{\bf Lemma}}\label{2: Lemma on repr of RGra in Cd(A)}
{\em For any (possibly, infinite-dimensional) graded vector space $A$
\Bi
\item[(i)]
the properad
$\ORGra_{d}$ admits a natural representation in the full Hochschild space
$\widehat{C}_{[d]}^\bu(A)$ of $A$,
$$
\rho_A: \ORGra_d \lon \cE nd_{\widehat{C}_{[d]}^\bu(A)},
$$
\item[(ii)] the operad $\DRTre_{d}$ admits a natural representation in the subspace
$C_{[d]}^\bu(A)\subset  \widehat{C}_{[d]}^\bu(A)$,
$$
\rho_A^0: \DRTre_d \lon \cE nd_{C_{[d]}^\bu(A)}.
$$
\Ei
}

\begin{proof} If $A$ is finite-dimensional, then $\widehat{C}_{[d]}^\bu(A)$ is identical to
the graded vector space
$$
 Cyc^\bu(W):= \prod_{k\geq 0} (\ot^k W)_{\Z_k}, \ \ \
  \text{where}\ W:= A[2-d]\oplus A^*[-1]
$$
The latter comes equipped with a { non-symmetric} pairing of degree $1-d$,
$$
< a_1\oplus b_1, a_2\oplus b_2>:= b_1(a_2), \ \ \
\forall\, a_1\oplus b_1,\ a_2\oplus b_2\in
  A[2-d] \oplus \Hom(A[1],\K)
$$
which can be used to construct a representation of $\ORGra_d$ in $Cyc^\bu(W)$ via 
exactly the same construction
as the one used in the proof of Theorem~4.2.2 in  \cite{MW1}.

\sip

If $A$ is infinite-dimensional, the  action $\rho_A(\Ga)$  of $\Ga\in \ORGra_{d}(m,n)$ on
$\ot^n \widehat{C}_{[d]}^\bu(A)$  can be described along the same lines but with the
replacement of the above pairing by the substitution of the output (say, $a_2\in A[2-d]$)
of one linear map from $\widehat{C}_{[d]}^\bu(A)$
into the input (say, $b_1\in A[1]$) of another linear map from  $\widehat{C}_{[d]}^\bu(A)$ 
and performing the composition of maps.

\sip

If a ribbon graph $\Ga\in \ORGra_{d}(m,n)$ has genus $\geq 1$, then  $\rho_A(\Ga)$
maps $\ot^n C_{[d]}^\bu(A)$ into $\ot^m \widehat{C}_{[d]}^\bu(A)$ in general,
not into $\ot^m {C}_{[d]}^\bu(A)$, so that we have to work with the {\em full}\, higher
Hochschild space of $A$ for consistency, and obtain claim (i).
\sip

If $\Ga$ has genus zero, then the values of the operation $\rho_A(\Ga)$ applied
to an element in $\ot^n C_{[d]}^\bu(A)$ belong to $\ot^m {C}_{[d]}^\bu(A)$ so that the
claim (ii) holds true.
\end{proof}

\sip

There is a morphism of properads
$$
i^{or}: \Lie_d \lon   \ORGra_d
$$
given on the Lie bracket generator of $\Lie_d$ by
$$
i^{or}:
\Ba{c}\begin{xy}
 <0mm,0.66mm>*{};<0mm,3mm>*{}**@{-},
 <0.39mm,-0.39mm>*{};<2.2mm,-2.2mm>*{}**@{-},
 <-0.35mm,-0.35mm>*{};<-2.2mm,-2.2mm>*{}**@{-},
 <0mm,0mm>*{\bu};<0mm,4.1mm>*{^{^{\bar{1}}}}**@{},
   <0.39mm,-0.39mm>*{};<2.9mm,-4mm>*{^{_2}}**@{},
   <-0.35mm,-0.35mm>*{};<-2.8mm,-4mm>*{^{_1}}**@{},
\end{xy}\Ea
 \lon \frac{1}{2}\left(
 \Ba{c}\resizebox{11mm}{!}{  \xy
 (3.5,4)*{^{\bar{1}}};
 (9,1)*+{_2}*\frm{o}="A";
 (0,1)*+{_1}*\frm{o}="B";
 \ar @{<-} "A";"B" <0pt>
\endxy} \Ea
 +(-1)^d
 \Ba{c}\resizebox{11mm}{!}{  \xy
 (3.5,4)*{^{\bar{1}}};
 (9,1)*+{_1}*\frm{o}="A";
 (0,1)*+{_2}*\frm{o}="B";
 \ar @{<-} "A";"B" <0pt>
\endxy} \Ea\right).
$$
Hence one can apply to $\ORGra_d$ Thomas Willwachers' twisting endofunctor \cite{Wi1} 
to obtain, in a full analogy to $\RGraphs_d$ above,  a {\em differential}\, graded properad
$$
\tw \ORGra_d=:\ORGraphs_d=\left\{ \ORGraphs_d(m,n), \delta   \right\}_{m\geq 1, n \geq 0}
$$
which is generated by ribbon graphs $\Ga$ with $m\geq 1$ labelled boundaries, $n\geq 0$ 
labelled vertices, and any
number of unlabelled vertices which are shown now in pictures in black color as $\bu$, 
and which are assigned the cohomological degree $d$. For example,
$$
\Ba{c} \resizebox{11mm}{!}{\xy
 (5,3)*{^{\bar{1}}};
 (10,1)*+{_1}*\frm{o}="B";
 (0,1)*{\bu}="A";
 \ar @{->} "A";"B" <0pt>
\endxy} \Ea
\in \ORGraphs_d(1,1), \ \
\Ba{c}\resizebox{17mm}{!}{
\mbox{$\xy
 (-2,3)*{^{\bar{1}}};
 (-2,-3)*{_{\bar{2}}};
 (-9,0)*{^{\bar{3}}};
 (-4,0)*{\bu}="C";
  (9,0)*+{_{_1}}*\frm{o}="1";
(-7,8)*+{_{_2}}*\frm{o}="2";
(-7,-8)*{\bu}="3";
 \ar @{->} "1";"C" <0pt>
  \ar @{->} "1";"2" <0pt>
   \ar @{<-} "1";"3" <0pt>
 \ar @{->} "2";"C" <0pt>
  \ar @{->} "3";"C" <0pt>
\endxy$}}
\Ea
\in \ORGraphs_d(3,2).
$$
The cohomological degree of $\Ga\in \ORGraphs_d$ is given by
$$
|\Ga|=(1-d) \# E(\Ga) + d\# V_\bu(\Ga),
$$
where $V_\bu(\Ga)$ stands for the set of unlabelled vertices.
The differential in $\ORGraphs_d(m,n)$ is given by the standard ``splitting of vertices" 
formula (cf.\ (\ref{2: delta in Tw(RGra)}))
\Beq\label{2: delta in Tw(RGra^or)}
\delta\Ga:=
\sum_{i=1}^m \left(\hspace{-2mm} \Ba{c}\resizebox{4mm}{!}{  \xy
 (0,7)*{\bu}="A";
 (0,0)*+{_1}*\frm{o}="B";
 \ar @{->} "A";"B" <0pt>
\endxy} \Ea
+(-1)^d\hspace{-2mm}
\Ba{c}\resizebox{4mm}{!}{  \xy
 (0,7)*{\bu}="A";
 (0,0)*+{_1}*\frm{o}="B";
 \ar @{<-} "A";"B" <0pt>
\endxy} \Ea
\right)
\  _1\hspace{-0.7mm}\circ_i \Ga\ \
- \ \ (-1)^{|\Ga|} 
\sum_{j=1}^n\Ga\  _j\hspace{-0.7mm}\circ_1
\left(\hspace{-2mm} \Ba{c}\resizebox{4mm}{!}{  \xy
 (0,7)*{\bu}="A";
 (0,0)*+{_1}*\frm{o}="B";
 \ar @{->} "A";"B" <0pt>
\endxy} \Ea
+(-1)^d\hspace{-2mm}
\Ba{c}\resizebox{4mm}{!}{  \xy
 (0,7)*{\bu}="A";
 (0,0)*+{_1}*\frm{o}="B";
 \ar @{<-} "A";"B" <0pt>
\endxy} \Ea
\right)
\  -\hspace{-2mm}\sum_{v\in V_\bu(\Ga)} \Ga\circ_v  \left(\xy
 (0,0)*{\bu}="a",
(5,0)*{\bu}="b",
\ar @{->} "a";"b" <0pt>
\endxy\right)
\Eeq

Lemma {\ref{2: Lemma on repr of RGra in Cd(A)}}(i) implies immediately the following result.

\subsection{Proposition} {\em Given any (possibly, curved) pre-CY structure $\pi$ on $A$, 
there is an associated action
of the dg properad $\ORGraphs_d$ on the full higher Hochschild complex 
$\left(\widehat{C}_{[d]}^\bu(A),\delta_\pi\right)$,
\Beq\label{2: action of OGrapgs on widehat{C}}
\rho^{full}_\pi: \ORGraphs_d \lon \cE nd_{\widehat{C}_{[d]}^\bu(A)}.
\Eeq
}

\sip

The above  action $\rho^{full}_\pi$ does not restrict to ${C}_{[d]}^\bu(A)$ in general, 
but this restriction works obviously fine for the sub-properad $(\DRTrees_d, \delta)$ 
of directed ribbon graphs of genus zero,
$$
\rho_\pi^0: \DRTrees_d \lon \cE nd_{{C}_{[d]}^\bu(A)},
$$
and hence induces a morphism of the associated cohomology operads,
$$
[\rho_\pi^0]: H^\bu(\DRTrees_d,\delta) \lon \cE nd_{H^\bu({C}_{[d]}^\bu(A),\delta_\pi)}.
$$

Recall that the collection of the homology groups
$$
\cG rav:=\{H_{\bu+1}(\cM_{0,n+1})\}_{n\geq 1}
$$
of the moduli space of genus zero algebraic curves with $n+1$ punctures has an operad structure
called the {\em gravity operad} \cite{Ge}.

\subsubsection{\bf Lemma}{\em The operad  $H^\bu(\DRTrees_d,\delta)$ is isomorphic 
to the degree shifted gravity operad},
$$
H^\bu(\DRTrees_d,\delta)\simeq  \cG rav\{d-2\}
$$
\begin{proof} Recall that  $\cR\cT rees_d$ stands for genus zero of the properad 
$\RGraphs_d$.  There is a morphism of dg operads
$$
t: \cR\cT rees_d \lon \DRTrees_d
$$
which sends every undirected edge of a ribbon tree from $\cR\cT rees_d$
 into a sum of two directed edges, e.g.
$$
\Ba{c} \resizebox{8mm}{!}{\xy
 (7,0)*+{_2}*\frm{o}="B";
 (0,0)*+{_1}*\frm{o}="A";
 \ar @{.} "A";"B" <0pt>
\endxy} \Ea \stackrel{t}{\lon}
\Ba{c} \resizebox{8mm}{!}{\xy
 (7,0)*+{_2}*\frm{o}="B";
 (0,0)*+{_1}*\frm{o}="A";
 \ar @{->} "A";"B" <0pt>
\endxy} \Ea
+(-1)^d
\Ba{c} \resizebox{8mm}{!}{\xy
 (7,0)*+{_2}*\frm{o}="B";
 (0,0)*+{_1}*\frm{o}="A";
 \ar @{<-} "A";"B" <0pt>
\endxy} \Ea
$$
We claim that this map is a quasi-isomorphism. Indeed, consider a filtration of both 
sides
by the number of vertices with valency $\geq 3$ and an associated series of morphisms 
at the level of the associated spectral sequences of both sides. The initial page of 
the left hand side of the map $t$ is a trivial
complex,   while the differential in the initial page of the right hand side creates
 only bivalent vertices. Using the argument identical to the one used in the proof of 
 Proposition K.1 in \cite{Wi1}
one concludes that the induced morphism of next pages of the spectral sequence is an 
isomorphism, i.e.
$H^\bu(\cR\cT rees_d)\simeq H^\bu(\DRTrees_d)$. It remains to cite the 
result of \cite{Wa}
which states an isomorphism of operads $
  H^\bu (\cR\cT rees_d) =  \cG rav\{d-2\}$.
\end{proof}

\subsubsection{\bf Corollary} {\em For any pre-CY algebra $(A,\pi)$ the associated
 higher Hochschild cohomology $H^\bu(C^\bu_{[d]}(A), \delta_\pi)$ is naturally a 
  $\cG rav\{d-2\}$-algebra.}

\sip

The Lie bracket part of $\cG rav\{d-2\}$-structure in $H^\bu(C^\bu_{[d]}(A), \delta_\pi)$ 
is controlled by the cohomology class
$$
\Ba{c} \resizebox{8mm}{!}{\xy
 (7,0)*+{_2}*\frm{o}="B";
 (0,0)*+{_1}*\frm{o}="A";
 \ar @{->} "A";"B" <0pt>
\endxy} \Ea
+(-1)^d
\Ba{c} \resizebox{8mm}{!}{\xy
 (7,0)*+{_2}*\frm{o}="B";
 (0,0)*+{_1}*\frm{o}="A";
 \ar @{<-} "A";"B" <0pt>
\endxy} \Ea \in H^\bu(\DRTrees_d(1,2))
$$
and hence coincides with the one induced from the standard Lie bracket in $C^\bu_{[d]}(A)$.
The next --- the trinary operation $[\ ,\ ,\ ]$ in $\cG rav\{d-2\}$-structure --- 
is controlled by
the following cohomology class
$$
\sum
\Ba{c}\resizebox{12mm}{!}{ \xy
(-7,0)*+{_1}*\frm{o}="0";
(7,0)*+{_2}*\frm{o}="1";
(0,12)*+{_3}*\frm{o}="2";
  (0,5)*{\bu}="s";
 \ar @{-} "s";"2" <0pt>
 \ar @{-} "s";"0" <0pt>
 \ar @{-} "s";"1" <0pt>
\endxy} \Ea \in H^\bu(\DRTrees_d(1,3))
$$
where the sum is taken over (skew)symmetrization of indices $1,2,3$ and making edges directed in 
all possible ways.

\subsection{Main definition-construction}\label{2: subsec on OORGraphs} Let $I^t$ be a 
dg ideal in 
$\ORGraphs_d$ generated by graphs having at least one black {\em target}, that is, a 
black (i.e.\  unlabelled) vertex with no
outgoing directed edges. Consider a dg sub-properad of the quotient properad,
$$
\OORGraphs_d\subset \ORGraphs_d/I^t,
$$
generated by oriented ribbon graphs with no outgoing edges at {\em labelled}\, vertices. 
It is easy to check that $\OORGraphs_d$ is a dg sub-properad indeed, i.e.\ it is 
closed under   both the  induced differential and the properadic compositions. 
Note that the generators  of $\OORGraphs_d$ 
must have at least one labelled white vertex, so that its decomposition into 
$\bS_m^{op}\times \bS_n$-modules
$$
\OORGraphs_d=\{ \OORGraphs_d(m,n)\}_{m\geq 1, n\geq1}
$$
has the ``input" integer parameter $n$ running from $1$ (not from zero as in the 
case of $\ORGraphs_d$).
Moreover, every generator of  $\OORGraphs_d$ except the unit $\resizebox{3.35mm}{!}{ \xy
 (0,1.5)*+{1}*\frm{o}="B";
\endxy}\in \ORGraphs_d(1,1)$ must have at least one black vertex.
For example,
$$
\Ba{c} \resizebox{11mm}{!}{\xy
 (5,3)*{^{\bar{1}}};
 (10,1)*+{_1}*\frm{o}="B";
 (0,1)*+{_2}*\frm{o}="A";
 \ar @{->} "A";"B" <0pt>
\endxy} \Ea
\not\in \OORGraphs_d,
\ \
\Ba{c}\resizebox{4mm}{!}{ \xy
 (0,8)*{\bu}="A";
 (0,0)*+{_1}*\frm{o}="B";
 \ar @{->} "A";"B" <0pt>
\endxy} \Ea
\in \OORGraphs_{d}(1,1),\ \
\Ba{c}\resizebox{15mm}{!}{
\mbox{$\xy
 (0,3)*{^{\bar{1}}};
 (-5,0)*{_{\bar{2}}};
 (-1,0)*{\bu}="C";
  (9,0)*+{_{_1}}*\frm{o}="1";
(-7,8)*+{_{_2}}*\frm{o}="2";
(-7,-8)*{\bu}="3";
 \ar @{<-} "1";"C" <0pt>
   \ar @{<-} "2";"3" <0pt>
 \ar @{<-} "2";"C" <0pt>
  \ar @{->} "3";"C" <0pt>
\endxy$}}
\Ea
\in \OORGraphs_d(2,2).
$$
The differential in $\OORGraphs_d$ is given on generators by
\Beq\label{2: delta in OORGra_d}
\delta\Ga:=
\sum_{i=1}^m \hspace{-1mm} \Ba{c}\resizebox{4mm}{!}{  \xy
 (0,7)*{\bu}="A";
 (0,0)*+{_1}*\frm{o}="B";
 \ar @{->} "A";"B" <0pt>
\endxy} \Ea
\  _1\hspace{-0.7mm}\circ_i \Ga\ \
- \ \ (-1)^{|\Ga|} 
\sum_{j=1}^n\Ga\  _j\hspace{-0.7mm}\circ_1
\hspace{-1mm} \Ba{c}\resizebox{4mm}{!}{  \xy
 (0,7)*{\bu}="A";
 (0,0)*+{_1}*\frm{o}="B";
 \ar @{->} "A";"B" <0pt>
\endxy} \Ea
\  -
\sum_{v\in V_\bu(\Ga)} \Ga\circ_v  \left(\xy
 (0,0)*{\bu}="a",
(5,0)*{\bu}="b",
\ar @{->} "a";"b" <0pt>
\endxy\right)
\Eeq
where in the last sum we set to zero every graph with at least one black target.

\mip

The following result is a cyclic version of the first half of Theorem 1 in \cite{KTV}.

\sip

\subsubsection{\bf Main Theorem} {\em For any degree $d$ pre-CY extension $\pi$ of a given 
$A_\infty$ structure
$\mu$ on a graded vector space $A$, there is an associated action of the dg properad
$\OORGraphs_d$ on the cyclic Hochschild complex $(Cyc^\bu(A,\K), \delta_\mu)$,
$$
\rho_\pi: \OORGraphs_d \lon \cE nd_{Cyc^\bu(A,\K)}.
$$
}

\begin{proof} Recall that under the assumptions of the theorem there is an action
(\ref{2: action of OGrapgs on widehat{C}})
the properad  $\ORGraphs_d$  on the full Hochschild complex
$(\widehat{C}_{[d]}^\bu(A), \delta_\pi)$ of the pre-CY algebra  $A$,
$$
\rho_\pi^{full}:  \ORGraphs_d  \lon \cE nd_{\widehat{C}_{[d]}^\bu(A)}
$$
obtained from Lemma~{\ref{2: Lemma on repr of RGra in Cd(A)}}(i) by the standard 
twisting procedure 
\cite{Wi1,Me2}. Under this representation,  the black vertices of the generators
 of $\ORGraphs_d$  
 get decorated with the given MC element element,
$$
\pi\in C_{[d]}^\bu(A) =\prod_{k\geq 1} C_{[d]}^{(k)}(A):= \prod_{k\geq 1}
\left( \bigoplus_{n_1,\ldots, n_k\geq 0}
\Hom\left(\bigotimes_{i=1}^{n_i}
(A[1])^{\ot n_i}, (A[2-d])^{\ot k}\right)_{\Z_k}\right),
$$
in such a way that every black vertex with precisely  $p$ outgoing edges is decorated 
with the summand
$\pi^{(\geq p)}$ of $\pi$ which belongs to the subspace
$
\prod_{k\geq p} C_{[d]}^{(k)}(A)\subset C_{[d]}^\bu(A)
$.
In particular, black targets in $\Ga$ (if any) {\em can}\, in principle contribute non-trivially 
into the representation $\rho_\pi^{full}$.

\sip

The acclaimed representation $\rho_\pi$  in the Theorem is built from $\rho_\pi^{full}$ via
 with the following crucial modification: a black vertex of a generator $\Ga\in \OORGraphs_d$
 with precisely $p$ outgoing edges is decorated now only with the summand $\pi^{(p)}$ of $\pi$
 which belongs to the subspace $C_{[d]}^{(p)}(A)\subset  \widehat{C}_{[d]}^\bu(A)$. As
 $\pi$ has no ``curvature" $\pi^{(0)}$ summand by its very definition, this modification of the 
 representation $\rho^{full}_\pi$ above factors through the ideal $I^t\subset \ORGraphs_d$ 
 generated by graphs having at least one target, that is, having at least one vertex with 
 $p=0$ outgoing edges.
 As white vertices of $\Ga$ have no out-going edges, the resulting modification $\rho_\pi$
 of $\rho^{full}_\pi$ is well-defined only on the summand $Cyc^\bu(A,\K)$ of 
 $\widehat{C}_{[d]}^\bu(A)$, 
 and that summand has to be equipped -- in accordance with the formula 
 (\ref{2: delta in OORGra_d})
  -- with the differential $\delta_\mu$ rather than with full differential $\delta_\pi$.
\end{proof}

\subsection{On induced $\Holie_d$ algebra structure }
Let $\Holie_d$ be a minimal resolution of $\Lie_d$.
It is a dg free operad whose (skew)symmetric generators,
\Beq\label{2: Lie_inf corolla}
\Ba{c}\resizebox{22mm}{!}{ \xy
(1,-5)*{\ldots},
(-13,-7)*{_1},
(-8,-7)*{_2},
(-3,-7)*{_3},
(7,-7)*{_{n-1}},
(13,-7)*{_n},
 (0,0)*{\bu}="a",
(0,5)*{}="0",
(-12,-5)*{}="b_1",
(-8,-5)*{}="b_2",
(-3,-5)*{}="b_3",
(8,-5)*{}="b_4",
(12,-5)*{}="b_5",
\ar @{-} "a";"0" <0pt>
\ar @{-} "a";"b_2" <0pt>
\ar @{-} "a";"b_3" <0pt>
\ar @{-} "a";"b_1" <0pt>
\ar @{-} "a";"b_4" <0pt>
\ar @{-} "a";"b_5" <0pt>
\endxy}\Ea
=(-1)^d
\Ba{c}\resizebox{23mm}{!}{\xy
(1,-6)*{\ldots},
(-13,-7)*{_{\sigma(1)}},
(-6.7,-7)*{_{\sigma(2)}},
(13,-7)*{_{\sigma(n)}},
 (0,0)*{\bu}="a",
(0,5)*{}="0",
(-12,-5)*{}="b_1",
(-8,-5)*{}="b_2",
(-3,-5)*{}="b_3",
(8,-5)*{}="b_4",
(12,-5)*{}="b_5",
\ar @{-} "a";"0" <0pt>
\ar @{-} "a";"b_2" <0pt>
\ar @{-} "a";"b_3" <0pt>
\ar @{-} "a";"b_1" <0pt>
\ar @{-} "a";"b_4" <0pt>
\ar @{-} "a";"b_5" <0pt>
\endxy}\Ea,
\ \ \ \forall \sigma\in \bS_n,\ n\geq2,
\Eeq
have degrees $1+d-nd$.
 The differential in $\Holie_d$ is given by
\Beq\label{2: d in Lie_infty}
\delta\hspace{-3mm}
\Ba{c}\resizebox{21mm}{!}{\xy
(1,-5)*{\ldots},
(-13,-7)*{_1},
(-8,-7)*{_2},
(-3,-7)*{_3},
(7,-7)*{_{n-1}},
(13,-7)*{_n},
 (0,0)*{\bu}="a",
(0,5)*{}="0",
(-12,-5)*{}="b_1",
(-8,-5)*{}="b_2",
(-3,-5)*{}="b_3",
(8,-5)*{}="b_4",
(12,-5)*{}="b_5",
\ar @{-} "a";"0" <0pt>
\ar @{-} "a";"b_2" <0pt>
\ar @{-} "a";"b_3" <0pt>
\ar @{-} "a";"b_1" <0pt>
\ar @{-} "a";"b_4" <0pt>
\ar @{-} "a";"b_5" <0pt>
\endxy}\Ea
=
\sum_{A\varsubsetneq [n]\atop
\# A\geq 2}\pm
%
%
\Ba{c}\resizebox{19mm}{!}{\begin{xy}
<10mm,0mm>*{\bu},
<10mm,0.8mm>*{};<10mm,5mm>*{}**@{-},
<0mm,-10mm>*{...},
<14mm,-5mm>*{\ldots},
<13mm,-7mm>*{\underbrace{\ \ \ \ \ \ \ \ \ \ \ \ \  }},
<14mm,-10mm>*{_{[n]\setminus A}};
<10.3mm,0.1mm>*{};<20mm,-5mm>*{}**@{-},
<9.7mm,-0.5mm>*{};<6mm,-5mm>*{}**@{-},
<9.9mm,-0.5mm>*{};<10mm,-5mm>*{}**@{-},
<9.6mm,0.1mm>*{};<0mm,-4.4mm>*{}**@{-},
<0mm,-5mm>*{\bu};
<-5mm,-10mm>*{}**@{-},
<-2.7mm,-10mm>*{}**@{-},
<2.7mm,-10mm>*{}**@{-},
<5mm,-10mm>*{}**@{-},
<0mm,-12mm>*{\underbrace{\ \ \ \ \ \ \ \ \ \ }},
<0mm,-15mm>*{_{A}}.
\end{xy}}
\Ea
\Eeq
If $d$ is even, all the signs above are equal to $-1$.

\subsubsection{\bf Proposition}\label{3: prop on F from holieb to RGraphs}  {\em  
There is a morphism of dg properads
$$
\Holie_d \lon \OORGraphs_{d+1}
$$
given on the generators by
$$
 \Ba{c}\resizebox{22mm}{!}{ \xy
(1,-5)*{\ldots},
(-13,-7)*{_1},
(-8,-7)*{_2},
(-3,-7)*{_3},
(7,-7)*{_{n-1}},
(13,-7)*{_n},
 (0,0)*{\bu}="a",
(0,5)*{}="0",
(-12,-5)*{}="b_1",
(-8,-5)*{}="b_2",
(-3,-5)*{}="b_3",
(8,-5)*{}="b_4",
(12,-5)*{}="b_5",
\ar @{-} "a";"0" <0pt>
\ar @{-} "a";"b_2" <0pt>
\ar @{-} "a";"b_3" <0pt>
\ar @{-} "a";"b_1" <0pt>
\ar @{-} "a";"b_4" <0pt>
\ar @{-} "a";"b_5" <0pt>
\endxy}\Ea
\lon
\sum_{\sigma\in \bS_n\atop i_k:=\sigma(k)}\frac{(-1)^{d\sigma}}{n}\
{\resizebox{21mm}{!}{
\xy
(-6,-9)*+{_{_{i_n}}}*\frm{o}="a1",
(-6,9)*+{_{_{i_2}}}*\frm{o}="a2",
(6,9)*+{_{_{i_3}}}*\frm{o}="a3",
(10,0)*+{_{_{i_4}}}*\frm{o}="a4",
(6,-9)*+{_{...}}*\frm{o}="a5",
(-10,0)*+{_{_{i_1}}}*\frm{o}="a6";
(0,0)*{\bu}="b",
\ar @{<-} "a1";"b" <0pt>
\ar @{<-} "a2";"b" <0pt>
\ar @{<-} "a3";"b" <0pt>
\ar @{<-} "a4";"b" <0pt>
\ar @{<-} "a5";"b" <0pt>
\ar @{<-} "a6";"b" <0pt>
\endxy}}\ \ \  \ \ \  \forall \ n\geq 2.
$$
}

Proof is a straightforward calculation; we omit the details. This result is
a ribbon properad analogue of a morphism of operads (6) studied in \S 3.5 of \cite{Wi2}.

\sip

The above Proposition implies Corollary {\ref{1: corollary on Holie}}.

\bip

\bip

{\Large
\section{\bf Cohomology of the properad $\OORGraphs_{d+1}$ and moduli spaces $\cM_{g, m+n}$}
}

\mip

\subsection{An extension of $\ORGraphs_{d+1}$}  Recall that the part $\ORGraphs_{d+1}(m,n)$
 of the
properad $\ORGraphs_{d+1}$ (see \S {\ref{2: subsec on RGra^or and RGraphs^or}})
is  generated by directed ribbon  graphs $\Ga$ with $m$ labelled boundaries, $n$ labelled white
vertices, and any number of black vertices satisfying the condition that the directed edges never
form a closed directed path in $\Ga$. The cohomological degree of a graph $\Ga$ is given by
$$
|\Ga|=(d+1)\# V_\bu(\Ga) - d \# E(\Ga)
$$
where $V_\bu(\Ga)$ is the set of black vertices of $\Ga$. Let  $\ORGraphs_{d,d+1}(m,n)$ be 
an extension of
$\ORGraphs_{d+1}(m,n)$ in which the generating ribbon graphs $\Ga$ are allowed to have any 
number of
{\em unlabelled white}\, vertices to which we assign the degree $d$ (as in the case of 
$\RGraphs_d$),
e.g.
$$
\Ba{c}\resizebox{13mm}{!}{ \xy
(-2.5,2.5)*{^{\bar{1}}},
(-7,4)*{^{\bar{2}}},
(3,5)*{^{\bar{3}}},
(-5,0)*+{_1}*\frm{o}="1";
 (5,0)*{\bu}="11";
 (5,10)*+{_2}*\frm{o}="2";
 (-5,10)*{\circ}="22";
 \ar @{->} "1";"22" <0pt>
 \ar @{->} "11";"22" <0pt>
 \ar @{->} "11";"22" <0pt>
  \ar @{->} "11";"1" <0pt>
   \ar @{->} "11";"2" <0pt>
  \ar @{->} "2";"22" <0pt>
\endxy} \Ea \in  \ORGraphs_{d,d+1}(3,2), \ \
\Ba{c}\resizebox{13mm}{!}{ \xy
(-2.5,2.5)*{^{\bar{1}}},
(-7,4)*{^{\bar{2}}},
(3,5)*{^{\bar{3}}},
(-5,0)*+{_1}*\frm{o}="1";
 (5,0)*{\bu}="11";
 (5,10)*{\circ}="2";
 (-5,10)*{\bu}="22";
 \ar @{<-} "1";"22" <0pt>
 \ar @{->} "11";"22" <0pt>
 \ar @{->} "11";"22" <0pt>
  \ar @{->} "11";"1" <0pt>
   \ar @{->} "11";"2" <0pt>
  \ar @{<-} "2";"22" <0pt>
\endxy} \Ea \in  \ORGraphs_{d,d+1}(3,1) .
$$
The cohomological degree of such a graph is given by
$$
|\Ga|=(d+1)\# V_\bu(\Ga) + d\# V_\circ(\Ga) - d \# E(\Ga)
$$
where $V_\circ(\Ga)$ is the set of  unlabelled white vertices of $\Ga$. The collection
$$
\ORGraphs_{d,d+1}=\{\ORGraphs_{d,d+1}(m,n)\}_{m\geq 1,n\geq 0}
$$
forms obviously a (non-differential) properad which contains  $\ORGraphs_{d+1}$ as a 
(non-differential) 
sub-properad generated by
graphs with no unlabelled white vertices. Moreover it contains a sub-properad
$$
\OORGraphs_{d,d+1}=\{\OORGraphs_{d,d+1}(m,n)\}_{m\geq 1,n\geq 0}
$$
generated by graphs with white vertices --- both labelled and unlabelled --- having no outgoing
edges. Note that, in contrast to the case of the properad  $\OORGraphs_{d+1}$ introduced in
\S {\ref{2: subsec on OORGraphs}}, the generators of  $\OORGraphs_{d,d+1}$ are allowed to have 
black vertices which are targets; for this reason the formula
(\ref{2: delta in Tw(RGra^or)}) does {\em not}\, define a differential in
$\OORGraphs_{d,d+1}$.
Hence to make $\OORGraphs_{d,d+1}$ into a {\em dg}\, properad the naive differential
(\ref{2: delta in Tw(RGra^or)}) has to be modified. The main point of this section is that such a 
modification is possible.

\sip

Consider a degree 1 derivation of the properad $\OORGraphs_{d,d+1}$
given by its action on the generating graphs as follows
$$
\delta_\bu \Ga:=- \left(\sum_{v\in V_{\circ}(\Ga)}\delta_{\circ\bu}^v\Ga + \Ga \circ_v 
\left( 
\sum_{k=1}^\infty \frac{1}{k}
{\resizebox{24mm}{!}{
\xy
(-4,-5)*{\circ}="a1",
(-4,5)*{\circ}="a2",
(4,5)*{\circ}="a3",
(4,-5)*{\circ}="a4",
(7,0)*{\circ}="a5",
(-7,0)*{\circ}="a6",
(0,0)*{\bu}="b",
(6,-9)*{\underbrace{\hspace{12mm}}_{k\ \text{edges}} \hspace{12mm}},
\ar @{<-} "a1";"b" <0pt>
\ar @{<-} "a2";"b" <0pt>
\ar @{<-} "a3";"b" <0pt>
\ar @{<-} "a4";"b" <0pt>
\ar @{<-} "a5";"b" <0pt>
\ar @{<-} "a6";"b" <0pt>
\endxy}}\right)
+
\sum_{w\in V_{\bu}(\Ga) } \Ga \circ_w \left( \resizebox{2.6mm}{!}{
\xy
 (0,4)*{\bullet}="a",
(0,-2)*{\bu}="b",
\ar @{->} "a";"b" <0pt>
\endxy}
\right)\right)
$$
where $\delta_{\circ\bu}^v \Ga$ is the ribbon graph obtained from $\Ga$ by making the 
white unlabelled 
vertex $v$ into a black unlabelled vertex, and $\circ_v (X)$ (resp.\ $\circ_w (X)$) means
the substitution of the unique boundary of the corresponding ribbon graph $X$ into the 
unlabelled white
vertex $v$ (resp., black vertex $w$) and performing the standard properadic-like composition.

\sip

It is a straightforward but tedious calculation (cf.\ \cite{Me3,MWW}) to check that 
$\delta_\bu^2=0$, 
i.e.\ that $\delta_\bu$ is a differential in $\OORGraphs_{d,d+1}$. We shall twist it 
using the fact 
that any
element
$$
\ga\in \OORGraphs_{d,d+1}(1,1)
$$
gives us a derivation $D_\ga$ of the properad $\OORGraphs_{d,d+1}$ given on a
generator $\Ga\in\OORGraphs_{d,d+1}(m,n)$  by
$$
D_\ga(\Ga)=\sum_{i=1}^m \ga\, {_1\circ_i}\Ga - (-1)^{|\Ga||\ga|} 
\sum_{j=1}^n \Ga{_j\circ_1}\ga.
$$

\subsubsection{\bf Lemma} {\em The degree 1 ribbon graph
$$
\ga:= \sum_{k=0}^\infty
{\resizebox{24mm}{!}{
\xy
(-4,-5)*{\circ}="a1",
(-4,5)*{\circ}="a2",
(4,5)*{\circ}="a3",
(4,-5)*{\circ}="a4",
(7,0)*{\circ}="a5",
(-7,0)*+{_1}*\frm{o}="a6";
(0,0)*{\bu}="b",
(6,-9)*{\underbrace{\hspace{12mm}}_{k+1\ \text{edges}} \hspace{12mm}},
\ar @{<-} "a1";"b" <0pt>
\ar @{<-} "a2";"b" <0pt>
\ar @{<-} "a3";"b" <0pt>
\ar @{<-} "a4";"b" <0pt>
\ar @{<-} "a5";"b" <0pt>
\ar @{<-} "a6";"b" <0pt>
\endxy}}
$$
satisfies the equation
$$
\delta_\bu \ga + \ga_{_1\circ_1}\ga =0.
$$
}

Proof is a straightforward calculation (cf.\ \cite{MWW}).

\mip

We conclude that the properad $\OORGraphs_{d,d+1}$ can be equipped with a differential
given on the generators by
$$
\delta \Ga:= \delta_\bu \Ga + D_\ga(\Ga).
$$

\subsection{Interpolation between $\RGraphs_d$ and $\OORGraphs_{d+1}$}
There is a diagram of short exact sequences of dg properads,
$$
\xymatrix{
&       &     0\ar[d]   &  &  &\\
 &       &     \OORGraphs_{d,d+1}^T\ar[d]   &  & &\\
 0\ar[r] & \OORGraphs_{d,d+1}^{ess} \ar[r] &
 \OORGraphs_{d,d+1} \ar[r]^{\pi_2}\ar[d]^{\pi_1} & \RGraphs_d\ar[r] & 0\\
  &       &     \OORGraphs_{d+1}\ar[d]   &  & &\\
   &       &    0   &  & &\\
}
$$
defined by the following two dg ideals of the dg properad $(\OORGraphs_{d,d+1}, \delta)$:
\Bi
\item[(i)] the dg ideal $\OORGraphs_{d,d+1}^T$ is generated by ribbon graphs having at least
one unlabelled vertex (of any colour) with no out-going edges, i.e.\ at least one unlabelled
target. As every unlabelled white vertex is a target, one can say the defining generators
of $\OORGraphs_{d,d+1}^T$ have at least one unlabelled white vertex or at least one black
target. It is obvious that the quotient complex $\OORGraphs_{d,d+1}/\OORGraphs_{d,d+1}^T$
is precisely the dg properad $\OORGraphs_{d+1}$ introduced in \S {\ref{2: subsec on OORGraphs}}.

\item[(ii)] Let us call a black vertex of a graph $\Ga\in \OORGraphs_{d,d+1}$
{\em inessential}\, if it is bivalent and has two outgoing edges, i.e.\ if it looks like this
(cf.\ \cite{MWW})
$$
    \text{inessential black vertex}:\ \ \resizebox{12mm}{!}{
\xy
(-7,1)*{}="a1";
(7,1)*{}="a2";
   (0,1)*{\bu}="b",
\ar @{<-} "a1";"b" <0pt>
\ar @{<-} "a2";"b" <0pt>
\endxy}
  $$
  All other black vertices in $\Ga$ (if any) are called {\em essential}. By definition,
  the dg ideal  $\OORGraphs_{d,d+1}^{ess}$ is generated by ribbon graphs having at least one
 essential black vertex. The quotient properad $\OORGraphs_{d,d+1}/\OORGraphs_{d,d+1}^{ess}$
 can be identified  with $\RGraphs_d$ after replacing each inessential black vertex connected
 to some  labelled or unlabelled white vertices  $u$ and $v$  with a dotted edge
  \[
\Ba{c}\resizebox{18mm}{!}{
\xy
%
 (0,0)*+{u}="a",
(8,0)*{\bullet}="b",
(16,0)*+{v}="c",
\ar @{<-} "a";"b" <0pt>
\ar @{->} "b";"c" <0pt>
\endxy}\Ea
\  \to \
\Ba{c}\resizebox{11mm}{!}{
\xy
%
 (0,0)*{u}="a",
(8,0)*{v}="b",
\ar @{.} "a";"b" <0pt>
\endxy}\Ea
\]
The possibility $u=v$ is not excluded.
\Ei

\subsection{Proposition} {\em The epimorphism
$$
\OORGraphs_{d,d+1} \stackrel{\pi_1}{\lon} \OORGraphs_{d+1}
$$
is a quasi-isomorphism of dg properads.}
\begin{proof} It is enough to show that
$
\ker \pi_1=\OORGraphs_{d,d+1}^T
$ is acyclic.
Let us consider a filtration of the latter by the total number of vertices.
The induced differential in the associated graded  $gr\ker \pi_1$
 acts by changing the color of unlabelled vertices from white to black,
 $\circ \rightarrow \bu$. By Maschke's Theorem, to prove the acyclicity of
 $gr\ker \pi_1$ it is enough to prove the acyclicity of its version
 $gr\ker \pi_1^{marked}$ in which the generating graphs $\Ga$
have all their edges and vertices marked (say totally ordered), but the type of vertices is not fixed.
Then the associated graded complex decomposes into a product of complexes
$$
gr\ker \pi_1^{marked} =\prod_{\Ga} C_\Ga
$$
parameterized by ribbon graphs $\Ga$ whose types of unlabelled vertices are not fixed, and
$$
C_\Ga:=\bigotimes_{v\in V(\Ga)} C_v
$$
where $C_v$ is
\Bi
\item[(i)] a one-dimensional trivial complex (generated by a black vertex) if the vertex $v$
has at least one outgoing solid edge in ${\Ga}$,
\item[(ii)] a two-dimensional
acyclic complex generated by one black and one white vertex if the vertex
$v$ has no outgoing solid edges in $\Ga$.
\Ei
 Any generating
graph $\Ga$  has least one vertex of  type (ii) as $\Ga$ has no closed paths of directed edges.
Hence each complex $C_\Ga$  has at least one acyclic tensor factor $C_v$ so that
it is itself acyclic implying the acyclicity of  $gr\ker \pi_1^{marked}$.
 The proposition is proven (cf.\ \S 4.1 in \cite{MWW}).
\end{proof}

\subsection{Proposition} {\em The epimorphism
$$
\OORGraphs_{d,d+1} \stackrel{\pi_2}{\lon} \RGraphs_{d}
$$
is a quasi-isomorphism of dg properads.}
\begin{proof}  It is enough to show that
$
\ker \pi_2=\OORGraphs_{d,d+1}^{ess}
$ is acyclic. Let us consider a filtration of $\OORGraphs_{d,d+1}^{ess}$ by the total
number of essential vertices and show that the associated graded complex
$gr_e \OORGraphs_{d,d+1}^{ess}$ is acyclic.

\sip
 Note that the differential in $gr_e \OORGraphs_{d,d+1}^{ess}$ can {\em not}\, decrease
 the number of unlabelled white vertices as the part $\delta_{\circ\bu}$ of the differential
 (which just changes the color)  increases the number of essential vertices by $+1$. Hence
 we can consider a filtration of $gr_e \OORGraphs_{d,d+1}^{ess}$ by the number of
 unlabelled white vertices and study the associated graded
 $$
 X:=gr_w gr_e \OORGraphs_{d,d+1}^{ess}.
 $$
 We claim that
 $$
 H^\bu(X)=0
 $$
 which implies the acyclicity of $gr_e\OORGraphs^{ess}_{d,d+1}$ and hence proves the
 proposition. Indeed, the induced differential $\delta_0$ in $X$ creates one new inessential
 vertex from an edge between two essential (black or white) vertices:
$$
\Ba{c}\resizebox{11mm}{!}{
\xy
%
 (0,0)*{\bullet}="a",
 (8,0)*+{_k}*\frm{o}="b",
\ar @{->} "a";"b" <0pt>
\endxy}\Ea
\  \to \
 \Ba{c}\resizebox{18mm}{!}{
\xy
%
 (0,0)*{\bullet}="a",
(8,0)*{\bullet}="b",
 (16,0)*+{_k}*\frm{o}="c"
\ar @{<-} "a";"b" <0pt>
\ar @{->} "b";"c" <0pt>
\endxy}\Ea, \ \
\Ba{c}\resizebox{11mm}{!}{
\xy
%
 (0,0)*{\bullet}="a",
(8,0)*{\circ}="b",
\ar @{->} "a";"b" <0pt>
\endxy}\Ea
\  \to \
 \Ba{c}\resizebox{18mm}{!}{
\xy
%
 (0,0)*{\bullet}="a",
(8,0)*{\bullet}="b",
 (16,0)*{\circ}="c"
\ar @{<-} "a";"b" <0pt>
\ar @{->} "b";"c" <0pt>
\endxy}\Ea
,
\quad
\quad
\quad
\Ba{c}\resizebox{11mm}{!}{
\xy
%
 (0,0)*{\bullet}="a",
(8,0)*{\bullet}="b",
\ar @{->} "a";"b" <0pt>
\endxy}\Ea
\  \to \
 \Ba{c}\resizebox{18mm}{!}{
\xy
%
 (0,0)*{\bullet}="a",
(8,0)*{\bullet}="b",
(16,0)*{\bullet}="c",
\ar @{<-} "a";"b" <0pt>
\ar @{->} "b";"c" <0pt>
\endxy}\Ea.
$$
Hence one can a consider a version $(X^{marked},\delta_0)$ of the complex $(X,\delta_0)$
which is a dg space generated by ribbon graphs whose all essential vertices are distinguished,
say totally ordered; in particular, we can stop distinguishing between labelled and
unlabelled white vertices in $X^{marked}$ and use one and the same symbol
$\resizebox{3.35mm}{!}{ \xy (0,1.5)*+{k}*\frm{o}="B";
\endxy}$, $k\in \N_{\geq 1}$, for their notation. Every generator $\Ga$ of $X^{marked}$
has at least one essential black vertex $\bu$  and at least one white vertex. Hence every
generator $\Ga$ has a white vertex
 $\resizebox{3.35mm}{!}{ \xy
 (0,1.5)*+{k}*\frm{o}="B";
\endxy}$  connected  to some $i$-labelled essential black vertex $\bu$ either directly by an edge,
$$
\xy
(0,2)*{^{i}},
%
 (0,0)*{\bullet}="a",
 (8,0)*+{_k}*\frm{o}="b";
\ar @{->} "a";"b" <0pt>
\endxy
$$
or via some non-essential vertex $\bu$ as in the following picture
$$
 \Ba{c}\resizebox{18mm}{!}{
\xy
(0,2)*{^{i}},
%
 (0,0)*{\bullet}="a",
(8,0)*{\bullet}="b",
 (16,0)*+{_k}*\frm{o}="c"
\ar @{<-} "a";"b" <0pt>
\ar @{->} "b";"c" <0pt>
\endxy}\Ea.
$$
We can assume without loss of generality that such a pair of essential vertices has
minimal possible labels, say $k=1$ and $i=2$, and then consider a filtration of $X^{marked}$
by the number of inessential vertices  which are connected to a essential vertex with
label $\geq 3$. The associated graded
$gr X^{marked}$ is the tensor product of a trivial complex
and the complex $C_{12}$ which controls the types of all possible ``edges" between
vertices 1 and 2. One has
$$
C_{12}=\bigoplus_{k\geq 1} \odot^k C
$$
where $C$ is a 2-dimensional complex generated by the following vectors,
$$
C=\text{span}\left\langle\Ba{c}\resizebox{11mm}{!}{
\xy
(0,2)*{^{2}},
%
 (0,0)*{\bullet}="a",
 (8,0)*+{_1}*\frm{o}="b";
\ar @{->} "a";"b" <0pt>
\endxy}\Ea
\ , \
\Ba{c}\resizebox{18mm}{!}{
\xy
(0,2)*{^{2}},
%
 (0,0)*{\bullet}="a",
(8,0)*{\bullet}="b",
 (16,0)*+{_1}*\frm{o}="c"
\ar @{<-} "a";"b" <0pt>
\ar @{->} "b";"c" <0pt>
\endxy}\Ea
 \right\rangle
$$
and equipped with the differential sending the first vector to the second one.
This complex is acyclic implying the acyclicity of $X^{matked}$ and hence of $X$.
The proof is completed (cf.\ Lemma 6.2.3 in \cite{Me3}).
\end{proof}

The above two Propositions imply main Theorem {\ref{1: Main theorem on H(OOGraphs)}} in the introduction.

\subsection{On the functorial meaning of the interpolating properad $\OORGraphs_{d,d+1}$}
Let consider again the twisted properad $\tw\RGra_{d+1}^{or}$ from \S {\ref{2: subsec on RGra^or and RGraphs^or}} equipped with the differential
$$
\delta_\bu (\Ga):= -\hspace{-2mm}\sum_{v\in V_\bu(\Ga)} \Ga\circ_v  \left(\xy
 (0,0)*{\bu}="a",
(5,0)*{\bu}="b",
\ar @{->} "a";"b" <0pt>
\endxy\right)
$$
rather than with the full twisted differential (\ref{2: delta in Tw(RGra)}). The dg properad
$(\tw\ORGra_{d+1},
\delta_\bu)$ contains a {\em dg}\, sub-properad $(\widehat{O}\RGraphs_{d+1}, \delta_\bu)$ generated by ribbon
graphs with no outgoing edges at then labelled (white) vertices. The difference between $(\widehat{O}\RGraphs_{d+1},
\delta_\bu)$ and $(\OORGraphs_{d+1}, \delta)$ is not only in their differentials, but also in the fact that the
unlabelled black vertices of the  generators of $\hat{O}\RGraphs_{d+1}$ are allowed to be targets.

\sip
Recall that there exists a dg operad $c\Holie_d$ governing {\em curved}\,  strongly homotopy Lie algebras; it is, by definition, a dg free operad generated by
(skew)symmetric $n$-corollas  (\ref{2: Lie_inf corolla}) of degrees $1+d-nd$ for any $n\geq 0$
(rather than for any  $n\geq 2$ as in the case of $\Holie_d$). The differential in $c\Holie_d$ is given formally by the same splitting formula (\ref{2: d in Lie_infty}) with the condition $\#A\geq 2$ replaced by $\#A \geq 0$.

\subsubsection{\bf Lemma} {\em There is a morphism of dg properads,
$$
(c\Holie_d, \delta) \lon  (\widehat{O}\RGraphs_{d+1},
\delta_\bu)
$$
given on generators by
$$
\Ba{c}\resizebox{22mm}{!}{ \xy
(1,-5)*{\ldots},
(-13,-7)*{_1},
(-8,-7)*{_2},
(-3,-7)*{_3},
(7,-7)*{_{n-1}},
(13,-7)*{_n},
 (0,0)*{\bu}="a",
(0,5)*{}="0",
(-12,-5)*{}="b_1",
(-8,-5)*{}="b_2",
(-3,-5)*{}="b_3",
(8,-5)*{}="b_4",
(12,-5)*{}="b_5",
\ar @{-} "a";"0" <0pt>
\ar @{-} "a";"b_2" <0pt>
\ar @{-} "a";"b_3" <0pt>
\ar @{-} "a";"b_1" <0pt>
\ar @{-} "a";"b_4" <0pt>
\ar @{-} "a";"b_5" <0pt>
\endxy}\Ea
\lon
\sum_{\sigma\in \bS_n\atop i_k:=\sigma(k)}\frac{(-1)^{d\sigma}}{n}\
{\resizebox{21mm}{!}{
\xy
(-6,-9)*+{_{_{i_n}}}*\frm{o}="a1",
(-6,9)*+{_{_{i_2}}}*\frm{o}="a2",
(6,9)*+{_{_{i_3}}}*\frm{o}="a3",
(10,0)*+{_{_{i_4}}}*\frm{o}="a4",
(6,-9)*+{_{...}}*\frm{o}="a5",
(-10,0)*+{_{_{i_1}}}*\frm{o}="a6";
(0,0)*{\bu}="b",
\ar @{<-} "a1";"b" <0pt>
\ar @{<-} "a2";"b" <0pt>
\ar @{<-} "a3";"b" <0pt>
\ar @{<-} "a4";"b" <0pt>
\ar @{<-} "a5";"b" <0pt>
\ar @{<-} "a6";"b" <0pt>
\endxy}} \ \ \  \ \ \  \forall  n\geq 0.
$$
}

Proof is a straightforward calculation (cf.\ Proposition  {\ref{3: prop on F from holieb to RGraphs}}).

\sip

Applying the standard twisting endofunctor \cite{Wi1, DSV, Me2} to the above map, one obtains precisely the above dg properad $\OORGraphs_{d,d+1}$ which plays the cental role in our proof of the Main Theorem {{\ref{1: Main theorem on H(OOGraphs)}}.

\bip


\bip

\def\cprime{$'$}

\end{document}